\documentclass[preprint,1p,times]{elsarticle}
\usepackage{amsmath, amssymb, amsthm}
\usepackage{graphicx}
\usepackage{color,soul}
\usepackage{float}
\usepackage{lineno}
\usepackage{pgf}
\usepackage{epstopdf}

\journal{ }
\begin{document}

\begin{frontmatter}
	
\title{Local and global dynamics of a fractional-order predator-prey system with habitat complexity and the corresponding discretized fractional-order system}
\author[rmv]{Shuvojit Mondal\corref{cor1}}
\ead{shuvojitmondal91@gmail.com}
\cortext[cor1]{Corresponding author}
\author[ju]{Milan Biswas\fnref{fn1}}
\author[cmbe]{Nandadulal Bairagi\fnref{fn2}}
\ead{nbairagi.math@jadavpuruniversity.in}
\address[rmv]{Department of Mathematics, Rabindra Mahavidyalaya\\ Hooghly-712401, India}
\address[ju]{Department of Mathematics, A.J.C. Bose College\\ Kolkata-700020, India}
\address[cmbe]{Centre for Mathematical Biology and Ecology \\ Department of Mathematics, Jadavpur University\\ Kolkata-700032, India.}

\begin{abstract}
This paper is focused on local and global stability of a fractional-order predator-prey model with habitat complexity constructed in the Caputo sense and corresponding discrete fractional-order system. Mathematical results like positivity and boundedness of the solutions in fractional-order model is presented. Conditions for local and global stability of different equilibrium points are proved. It is shown that there may exist fractional-order-dependent instability through Hopf bifurcation for both fractional-order and corresponding discrete systems. Dynamics of the discrete fractional-order model is more complex and depends on both step length and fractional-order. It shows Hopf bifurcation, flip bifurcation and more complex dynamics with respect to the step size. Several examples are presented to substantiate the analytical results.
\end{abstract}


\end{frontmatter}

\section{Introduction}
\label{intro}
Fractional calculus is the area of mathematics where derivatives and integrals can be extended to an arbitrary order. There are different approaches to study the dynamical behaviors of population models, e.g. ordinary differential equations (ODE), partial differential equations (PDE), difference equations (DE), fractional-order differential equations (FDE) etc. The first three techniques are being extensively used for a long time. However, the fractional-order differential equations have gained considerable importance only in the recent past due to their ability of providing an exact or approximate description of different nonlinear phenomena. The main advantage of fractional-order system is that they allow greater degrees of freedom than an integer order system \cite{DasGupta11}. FDE are preferably used since they are naturally related to systems with memory which exists in most biological phenomena \cite{Ahmed07}. Moreover, FDE has close relations to fractals which has wide applications in mathematical biology. Recently, some authors have investigated the importance of fractional-order differential equations in several biological systems, e.g. ecological system with delay \cite{Rihan15,Zhang19}, control based epidemiological system \cite{Basir16}, ecological system with diffusion \cite{Zhang17} etc. It also has applications in other fields of science and engineering \cite{JV07,TET09,GG03,MM07,D11}. Some recent studies discuss about the approximate solution of nonlinear fractional-order differential population models \cite{DasGupta11,MET17} and some others study the qualitative behavior of nonlinear interactions of biological systems \cite{LET16,Vargas15,HET15,Bairagi17,Ghaziania16}.  However, existence and proof of Hopf bifurcation, that causes oscillations in population densities due to fractional-order, is rare in the fractional-order population model. We address this issue along with others in a fractional-order predator-prey model considered in Caputo sense. Lot of discrete models on biological systems have been proposed and analyzed. However, discretization of a fractional-order population model is rare. Elsadany and Matouk \cite{Matouk15} recently studied a fractional-order Lotka-Volterra predator prey model with its discretization. They showed complex dynamics even in a simpler prey-predator model. In the second phase of this paper, we construct the discrete version of the continuous fractional-order system and reveal its dynamics.\\

All most every habitat, whether it is aquatic or terrestrial, contains some kind of complexity. For example, sea grass, aquatic weeds, salt marshes, littoral zone vegetation, mangroves, coral reefs etc. make aquatic habitat complex. Both field and laboratory experiments confirm that habitat complexity increases persistency of interacting species \cite{August83,Beukers97,Canion09,Ellner01,Frederick06,Johnson03}. A general hypothesis is that habitat complexity reduces predation rates by decreasing predator-prey interaction and thereby increases population persistency. A Rosenzweig-MacArthur predator-prey model \cite{Rosenzweig63} that incorporates the effect of habitat complexity can be represented by the following coupled nonlinear system:
\begin{eqnarray}\label{Ecological_model}
\frac{dx}{dt} & = & rx\bigg(1-\frac{x}{K}\bigg) - \frac{\alpha (1-c) x y}{1+\alpha (1-c)h x}, x(0) > 0, \\
\frac{dy}{dt}& = & \frac{\theta \alpha (1-c) x y}{1+\alpha (1-c)h x} - dy, y(0) > 0. \nonumber
\end{eqnarray}
This model says that the prey population $x$ grows logistically with intrinsic growth rate $r$ to its carrying capacity $K$. Predator $y$ captures the prey at a maximum rate $\alpha$ in absence of any habitat complexity $(c=0)$. In presence of complexity, predation rate decreases to $\alpha(1-c)$, where the dimensionless parameter $c$ is called the degree or strength of complexity. The value of $c$ ranges from $0$ to $1$. In particular, $c = 0.4$ implies that predation rate decreases by $40 \%$ due to habitat complexity. If $c=0$, i.e. if there is no habitat complexity then the system (\ref{Ecological_model}) reduces to well known Rosenzweig-MacArthur model \cite{Rosenzweig63}. However, if $c=1$ then $y\rightarrow 0$ as $t\rightarrow \infty$ and the prey population grows logistically to its maximum value $K$. The parameter $\theta$ $(0<\theta<1)$ is the conversion efficiency, measuring the number of newly born predators for each captured prey and $d$ is the death rate of predator. All parameters are assumed to be positive. For construction and  more explanation of the model, readers are referred to \cite{Bairagi11}.\\

Considering the fractional derivatives in the sense of Caputo, we have the following fractional-order model corresponding to the integer order model (\ref{Ecological_model}):
\begin{eqnarray}\label{Ecological fractional order model}
^{c}_{0} D^{m}_{t}x & = & rx\bigg(1-\frac{x}{K}\bigg) - \frac{\alpha (1-c) x y}{1+\alpha (1-c)h x}, \\
^{c}_{0} D^{m}_{t}y & =  & \frac{\theta \alpha (1-c) x y}{1+\alpha (1-c)h x} - dy, \nonumber
\end{eqnarray}
where $^{c}_{0} D^{m}_{t}$ is the Caputo fractional derivative with fractional-order $m$ $(0< m \leq1)$. The main advantage of Caputo's approach is that the initial conditions for the fractional differential equations with Caputo derivatives takes the similar form as for integer-order differential equations \cite{CuiYang14,Podlubny99}, and thus takes the advantage of defining integer order initial conditions for fractional-order differential equations. We analyze system (\ref{Ecological fractional order model}) with the initial conditions $x(0)>0,~ y(0)>0.$ In this paper, we prove different mathematical results, like existence, non-negativity and boundedness of the solutions of fractional-order system (\ref{Ecological fractional order model}). We establish conditions for local and global stabilities of different equilibrium points. It is shown that the interior equilibrium may switch its stability through Hopf bifurcation for some critical value of the fractional-order when the degree of complexity is low. A discrete system generally produces more complex dynamics than its continuous counterpart \cite{Matouk15}. Here we construct a discrete version of the fractional-order prey-predator model (\ref{Ecological fractional order model}). We prove local stability of different fixed points of the discrete system along with the existence conditions of Hopf and flip bifurcations. Numerical examples are presented for both systems in support of the analytical results. He and Lai (2011) has discretized a continuous type predator-prey model by Euler method. Using center manifold theorem, it is shown that the system undergoes flip and Neimark-Sacker bifurcations. Period doubling bifurcation leading to chaos was also shown through numerical simulations. However, they have not studied the dynamics of fractional-order discrete system. Abdelaziz et al. (2018) transformed an integer order SI-type epidemic model to a fraction-order discrete epidemic model and analyzed it to show flip and Neimark-Sacker bifurcations. But did not analyze the qualitative behavior of the fractional-order system. Here we study both the fractional-order and discretized fractional-order predator-prey systems. We compare the qualitative behavior of integer order system with the fractional-order and discretized fractional-order systems.\\

The rest of the paper is organized as follows. The next section contains well-posedness, existence and uniqueness of the solutions of the fractional-order system. Qualitative behavior of different equilibrium points are also presented here. Section 3 deals with stability and hopf bifurcation of fractional-order discrete system. Different examples are presented to illustrate the observed dynamics in Section 4. The paper ends with a brief discussion in Section 5.
\section{Well-posedness}
\label{sec:1}
\subsection{\textbf{Nonnegativity and boundedness}}
Considering the biological significance of the model, we are only interested in solutions that are nonnegative and bounded in the region $\Re^2{}_{+} = \{z\in\Re^{2}|z\geq0 \}$ and $z(t) = (x(t),y(t))^{T}$. To prove the nonnegativity and  uniform boundedness of our system, we shall use the following results. \\

\noindent\textbf{Lemma 2.1} \cite{Odibat07} \textit{Suppose that $f(t)\in C[a,b]$ and $D^{m}_{a}f(t)\in C(a,b]$ with $0<m\leq1$. The Generalized Mean Value Theorem states that
	\begin{equation}\nonumber
	f(t) = f(a) + \frac{1}{\Gamma(m)}(D^{m}_{a}f)(\xi). (t-a)^m,
	\end{equation}
	where $a\leq\xi\leq t$, $\forall t \in (a,b]$.}\\

\noindent From this lemma, one can easily prove the following result.\\

\noindent\textbf{Corollary 2.1} \cite{LET16,Odibat07} \textit{ Suppose $f(t)\in C[a,b]$ and $^{c}_{t_0} D^{m}_{t}f(t)\in C[a,b]$, $0<m \leq 1$. If $^{c}_{t_0} D^{m}_{t}f(t)\geq 0, \forall t\in (a,b)$ then $f(t)$ is a non decreasing function for each $t\in[a,b]$ and if $^{c}_{t_0} D^{m}_{t}f(t)\leq 0, \forall t\in (a,b)$ then $f(t)$ is a non-increasing function for each $t\in[a,b]$.}\\

\noindent\textbf{Lemma 2.2} \cite{LET16} \textit{ Let $u(t)$ be a continuous function on $[t_{0},\infty)$ and satisfying
	$$^{c}_{t_0} D^{m}_{t}u(t) \leq -\lambda u(t) + \mu,$$
	$$u(t_{0})                         =  u_{t_{0}}, $$
	where $0< m\leq1$, $(\lambda,\mu)\in\Re^{2}$, $\lambda\neq0$ and $t_{0}\geq0$ is the initial time. Then its solution has the form
	\begin{equation}\nonumber
	u(t) \leq \bigg(u_{t_{0}}- \frac{\mu}{\lambda}\bigg)E_{m}[-\lambda(t-t_{0})^{m}] + \frac{\mu}{\lambda}.
	\end{equation}}
\noindent\textbf{Theorem 2.1} \textit{ All solutions of system (\ref{Ecological fractional order model}) which start in $\Re^{2}_{+}$ are nonnegative and uniformly bounded.}\\

\noindent \textbf{Proof} First we show that the solutions $x(t) \in \Re^{2}_{+} $ are nonnegative if it start with positive initial values. If not, then there exists a $t_{1}> 0$ such that
\begin{eqnarray} \label{Function}\nonumber
&x(t)&> 0,~~0\leq t< t_{1},\\
&x(t)& = 0, ~~ t = t_1, \\
&x(t_{1}^+)& < 0.\nonumber
\end{eqnarray}
Using (\ref{Function}) in the first equation of (\ref{Ecological fractional order model}), we have
\begin{equation}
^{c}_{0} D^{m}_{t} x(t)|_{t = t_{1}} = 0.
\end{equation}
According to Corollary $2.1$, we have $x(t_{1}^{+}) = 0$, which contradicts the fact $x(t_{1}^{+}) < 0$. Therefore, we have $x(t)\geq0,~\forall ~t\geq0$. Using similar arguments, we can prove $y(t)\geq0, \forall t\geq0$.\\ Next we show that all solutions of system (\ref{Ecological fractional order model}) which initiate in $\Re^{2}_{+}$ are uniformly bounded. Define a function
\begin{equation}
V(t) = x + \frac{1}{\theta}y,
\end{equation}
Taking fractional time derivative, we have
\begin{equation}\nonumber
^{c}_{0} D^{m}_{t} V(t) =~ ^{c}_{0} D^{m}_{t} x(t)+~ ^{c}_{0} D^{m}_{t} \frac{1}{\theta} y(t) ~=~ rx\bigg(1-\frac{x}{K}\bigg) - \frac{d}{\theta} y.
\end{equation}
Now, for each $\eta>0$, we have
\begin{align}\label{Boundedness}
^{c}_{0} D^{m}_{t} V(t) + \eta V(t) =  & rx\bigg(1-\frac{x}{K}\bigg) - \frac{d}{\theta} y + \eta x + \frac{\eta}{\theta} y \nonumber \\
=  &-\frac{r}{K}x^{2} + (r + \eta)x + (\eta - d)\frac{1}{\theta}y \nonumber\\
\leq & \frac{K}{4r}(r + \eta)^{2} + (\eta - d)\frac{1}{\theta}y.
\end{align}
If we take $\eta < d$ then right hand side of (\ref{Boundedness}) is bounded in $\Re^{2}_{+}$ and there exist a constant $l>0$ (say) such that
\begin{equation}
^{c}_{0} D^{m}_{t} V(t) + \eta V(t) \leq l,
\end{equation}
where $l  = \frac{K}{4r}(r + \eta)^{2}$.\\
Applying Lemma $2.2$, we then have
\begin{eqnarray}\nonumber
V(t) &\leq & (V(0) - \frac{l}{\eta})E_{m}[-\eta t^{m}] + \frac{l}{\eta} \nonumber \\
&\leq &  V(0) E_{m}[-\eta t^{m}] + \frac{l}{\eta}(1 - E_{m}[-\eta t^{m}]).
\end{eqnarray}
For $t\rightarrow \infty$, we thus have $V(t)\rightarrow \frac{l}{\eta}$. Therefore, $0<V(t)\leq \frac{l}{\eta}$. Hence all solutions of the system (\ref{Ecological fractional order model}) that starts from $\Re^{2}_{+}$ are confined in the region $B = \{(x,y)\in\Re^{2}_{+}| 0<V(t)\leq \frac{l}{\eta} + \epsilon$, for any $\epsilon>0, 0<\eta<d, l = \frac{K}{4r}(r + \eta)^{2}\}$. Hence the theorem.
\subsection{\textbf{Existence and uniqueness}}
Here we  study the existence and uniqueness of the solution of our system (\ref{Ecological fractional order model}). We have the following Lemma due to Li et al \cite{LiChen10}.

\noindent\textbf{Lemma 2.3} \textit{Consider the system
	\begin{equation}\nonumber
	^{c}_{t_{0}} D^{m}_{t} x(t) = f(t,x), t>t_{0}
	\end{equation}
	with initial condition $x_{t_{0}}$, where $0<m\leq1$, $f:[t_{0},\infty)\times\Omega\rightarrow\Re^{n}$, $\Omega\in\Re^{n}$. If $f(t,x)$ satisfies the locally Lipschitz condition with respect to $x$ then there exists a unique solution of the above system on $[t_{0},\infty)\times\Omega$.}

We study the existence and uniqueness of the solution of system (\ref{Ecological fractional order model}) in the region $\Omega \times[0, T]$, where $\Omega = \{(x,y)\in\Re^{2}|~ max\{|x|, |y|\} \leq M\}$, $T<\infty$ and $M$ is large. Denote $X = (x,y)$, $\bar{X} = (\bar{x}, \bar{y})$. Consider a mapping $H: \Omega\rightarrow\Re^{2}$ such that $H(X) = (H_{1}(X), H_{2}(X))$, where
\begin{equation}\label{Existence}
H_{1}(X)  =  rx\bigg(1-\frac{x}{K}\bigg) - \frac{\alpha (1-c) x y}{1+\alpha (1-c)h x},~H_{2}(X)  =  \frac{\theta \alpha (1-c) x y}{1+\alpha (1-c)h x} - dy.
\end{equation}
For any $X, \bar{X} \in \Omega$, it follows from (\ref{Existence}) that
\begin{eqnarray}\nonumber
\begin{split}
\parallel H(X) - H(\bar{X})\parallel = \mid H_{1}(X)-H_{1}(\bar{X})\mid +\mid H_{2}(X)-H_{2}(\bar{X})\mid  \nonumber \\
= \mid rx\bigg(1-\frac{x}{K}\bigg) - \frac{\alpha (1-c) x y}{1+\alpha (1-c)h x} - r\bar{x}\bigg(1-\frac{\bar{x}}{K}\bigg) + \frac{\alpha (1-c) \bar{x}\bar{y}}{1+\alpha (1-c)h\bar{x}}\mid  \\
+\mid \frac{\theta \alpha (1-c) x y}{1+\alpha (1-c)h x} - dy - \frac{\theta \alpha (1-c) \bar{x} \bar{y}}{1+\alpha (1-c)h \bar{x}} + d\bar{y}\mid  \\
= \mid r(x-\bar{x}) - \frac{r}{K}(x^{2} -\bar{x}^{2}) -\alpha (1-c) \bigg(\frac{xy}{1+ \alpha (1-c)hx} - \frac{\bar{x}\bar{y}}{1 + \alpha (1-c)h\bar{x}}\bigg)\mid  \\
+\mid \theta \alpha (1-c) \bigg(\frac{xy}{1+ \alpha (1-c)hx} - \frac{\bar{x}\bar{y}}{1 + \alpha (1-c)h\bar{x}}\bigg) - d(y - \bar{y})\mid  \\
\leq r\mid x-\bar{x}\mid + \frac{r}{K}\mid x^{2} -\bar{x}^{2}\mid \\
+\alpha (1-c)(1+\theta) \mid\bigg(\frac{xy}{1+ \alpha (1-c)hx} - \frac{\bar{x}\bar{y}}{1 + \alpha (1-c)h\bar{x}}\bigg)\mid + d\mid y - \bar{y}\mid  \\
\leq   r\mid x-\bar{x}\mid + \frac{2rM}{K}\mid x - \bar{x}\mid + \alpha(1-c)(1+\theta) \mid xy - \bar{x}\bar{y}\mid \\
+ \alpha^2 (1-c)^2 (1+\theta)hM^2 \mid y - \bar{y}\mid + d\mid y - \bar{y}\mid\\
\leq   \bigg(r +\frac{2rM}{K} + \alpha (1-c)(1+\theta)M \bigg)\mid x - \bar{x}\mid \\
+ \bigg(\alpha (1-c)(1+\theta)M + d + \alpha^2 (1-c)^2 (1+\theta)hM^2 \bigg)\mid y-\bar{y}\mid \\
\leq  L\parallel (x,y) - (\bar{x}, \bar{y})\parallel \\
\leq  L\parallel X - \bar{X} \parallel,
\end{split}
\end{eqnarray}
where $L = max \{r +\frac{2rM}{K} + \alpha (1-c)(1+\theta)M, \alpha (1-c)(1+\theta)M + d + \alpha^2 (1-c)^2 (1+\theta)hM^2\}$. Thus $H(X)$ satisfies Lipschitz condition with respect to $X$ and following Lemma $2.3$, there exists a unique solution $X(t)$ of system (\ref{Ecological fractional order model}) with initial condition $X(0) = (x(0), y(0))$.
\subsection{\textbf{Stability of equilibrium points}}
We have the following stability result on fractional-order differential equations.\\

\noindent\textbf{Theorem 2.2} \cite{Petras11} \textit{Consider the following fractional-order system
	\begin{equation}\label{Stability condition} \nonumber
	^{c}_{0} D^{m}_{t} x(t) = f(x), x(0)= x_{0}
	\end{equation}
	with $0<m\leq1, x\in \Re^{n}$ and $f: \Re^{n}\rightarrow \Re^{n}$. The equilibrium points of the  above system are calculated by solving the equation $f(x) = 0$. These equilibrium points are locally asymptotically stable if all eigenvalues $\lambda_{i}$ of the jacobian matrix $J = \frac{\partial f}{\partial x}$ evaluated at the equilibrium points satisfy $$\mid arg(\lambda_{i})\mid >\frac{m \pi}{2}, i = 1,2,.....,n.$$}
For any quadratic polynomial $\phi(x) = x^2 + a_1 x + a_2$, the discriminant $D(\phi)$ of the polynomial $\phi$ is given by
\[
\mathbf{D(\phi)} = - \begin{vmatrix}
1 & a_{1} & a_{2} \\ 2 & a_{1} & 0 \\ 0 & 2& a_1
\end{vmatrix} = a_1^2 - 4a_2.
\]
The generalized Routh-Hurwitz stability conditions for fractional-order systems are then given by the following proposition \cite{Ahmed07,Ahmed05,Ahmed06}.\\

\noindent\textbf{Proposition 2.1} \textit{\begin{itemize} \item[(i)] If $D(\phi) \geq 0$, $ a_{1}>0$ and $a_{2}>0$, then the equilibrium $E^{*}$ is locally asymptotically stable for $0< m \leq 1$.
		\item[(ii)] If $D(\phi) < 0$, $a_{1} < 0$ and $\mid \tan^{-1}(\frac{\sqrt{4a_2 - a_1^{2}}}{a_1}) \mid > \frac{m \pi}{2}$, $0< m < 1$, then the equilibrium $E^{*}$ is locally asymptotically stable.
\end{itemize}}

\noindent The system (\ref{Ecological fractional order model}) has three equilibrium points: $(i)$ $E_{0} = (0,0)$ as the trivial equilibrium, $(ii)$ $E_{1} = (K,0)$ as the predator-free equilibrium and $(iii)$ $E^{*} = (x^{*},y^{*})$ as the interior equilibrium, where
\begin{eqnarray}\label{Equilibrium relation}
x^{*} & = & \frac{d}{\alpha (1-c)(\theta - h d)}, ~y^{*} = \frac{r(K-x^*)\{1 + \alpha h(1-c)x^{*}\}}{\alpha K(1-c)}.
\end{eqnarray}
Note that the equilibria $E_{0}$ and $E_{1}$ always exist. The interior equilibrium $E^{*}$ exists if $0 < c < c_{1}$ and $\theta > \theta_{1}$, where $c_{1} = 1 - \frac{d}{\alpha K (\theta - h d)}$, $\theta_{1} = hd + \frac{d}{\alpha K}$.\\

\noindent\textbf{Theorem 2.3} \textit{(a) The trivial equilibrium point $E_{0}$ is a saddle point. (b) The predator-free equilibrium point $E_{1}$ is locally asymptotically stable if $c>c_{1}$, $\theta > \theta_{1}$ and a saddle if $c<c_{1}.$ } \\

\noindent\textbf{Proof}  The proof of part (a) is straightforward and omitted. The Jacobian matrix corresponding to $E_{1}$ is given by

%

\[
\mathbf{J(E_{1})} = \begin{pmatrix}
-r & -\frac{\alpha K(1-c)}{1+\alpha K(1-c)h} \\ 0 & \frac{\theta \alpha K(1-c)}{1+\alpha K(1-c)h}-d \end{pmatrix}.\]
The corresponding eigenvalues are $\xi_{1} = -r ~(<0)$, $\xi_{2} = \frac{\theta \alpha K(1-c)}{1+\alpha K(1-c)h} - d$. If $c<c_{1}$, then $\xi_{2}>0$ and $\mid arg(\xi_{2})\mid = 0$. In this case, $E_{1} = (K, 0)$ is a saddle point. \\

If $c>c_{1}$ and $\theta > \theta_{1}$ then $\xi_{2}<0$. Consequently, $\mid arg(\xi_{i})\mid  = \pi > \frac{m \pi}{2}, \forall m\in (0,1]$, $i = 1, 2,$ and the equilibrium $E_{1} = (K, 0)$ is  locally asymptotically stable. In other words, when the degree of complexity is high and the conversion efficiency of predator exceeds some lower threshold value, then the predator-free equilibrium becomes locally asymptotically stable. \\

To prove the global stability of $E_{1}$, we use the following Lemma. \\

\noindent\textbf{Lemma 2.4} \cite{Vargas15} \textit{Let $x(t)\in \Re_{+}$ be a continuous and derivable function. Then for any time instant $t>t_{0}$
	\begin{equation}\nonumber
	^{c}_{t_{0}} D^{m}_{t}\bigg[x(t) - x^{*} - x^{*}ln\frac{x(t)}{x^{*}}\bigg] \leq \bigg(1-\frac{x^{*}}{x(t)}\bigg)~~{^{c}_{t_{0}}} D^{m}_{t}x(t), x^{*}\in \Re_{+}, \forall m\in(0,1].
	\end{equation}}

\noindent\textbf{Theorem 2.4} \textit{The predator-free equilibrium $E_{1}$ is globally asymptotically stable for any $m\in (0, 1]$ if $c > c_{1}$, $\theta > \theta_{1}$, where $c_{1} = 1 - \frac{d}{\alpha K (\theta - h d)}$, $\theta_{1} = hd + \frac{d}{\alpha K}$.}\\

\noindent\textbf{Proof}
Consider the Lyapunov function
\begin{equation}
V(x,y) = \bigg(x - K - K ln\frac{x}{K}\bigg) + \frac{y}{\theta}.
\end{equation}
Here $V(x,y) > 0$ for all values of $x(t), y(t)> 0$ and $V = 0$ only at $E_{1} = (K,0)$. Calculating the $mth$ order fractional derivative of $V(x,y)$ along the solution of (\ref{Ecological fractional order model}) and using Lemma $2.4$ when $t_{0} = 0$, we have
\begin{eqnarray}\label{Globality_1}\nonumber
\begin{split}
^c_{0}D^{m}_{t}V(x,y) \leq& \frac{(x - K)}{x} {^c_{0}}D^{m}_{t}x(t) + \frac{1}{\theta}{^c_{0}}D^{m}_{t}y(t)\\
=& (x - K)\bigg[r(1-\frac{x}{K}) - \frac{\alpha (1-c) y}{1+\alpha (1-c)h x}\bigg] +\frac{\alpha (1-c) x y}{1+\alpha (1-c)h x} - \frac{dy}{\theta} \\
=& (x - K)\bigg[-\frac{r}{K}(x-K) - \frac{\alpha (1-c) y}{1+\alpha (1-c)h x}\bigg] + \frac{\alpha (1-c) x y}{1+\alpha (1-c)h x} - \frac{dy}{\theta}\\
=& -\frac{r}{K}(x - K)^{2} + \frac{\alpha K(1-c) y}{1+\alpha (1-c)h x} - \frac{dy}{\theta}\\
\leq & -\frac{r}{K}(x - K)^{2} + \bigg[\alpha K(1-c) - \frac{d}{\theta}\bigg]y.\\
\end{split}
\end{eqnarray}
One can note that $^c_{0}D^{m}_{t}V(x,y) \leq0, \forall (x,y)\in R^{2}_{+}$ if $\alpha K(1-c) - \frac{d}{\theta} < 0$, i.e., if $\frac{d}{\theta} > \alpha K(1-c) > \frac{\alpha K(1-c)}{1+\alpha (1-c)h K}$. This implies $^c_{0}D^{m}_{t}V(x,y) \leq0, \forall (x,y)\in R^{2}_{+}$ if $c > c_{1}$, $\theta > \theta_{1}$ and $^c_{0}D^{m}_{t}V(x,y) = 0$ at $E_{1}$. Therefore, the only invariant set on which $^c_{0}D^{m}_{t}V(x,y) = 0$ is the singleton $\{E_{1}\}$. Then by Lemma $4.6$ in \cite{HET15}, it follows that the predator-free equilibrium $E_{1}$ is globally asymptotically stable if $c > c_{1}$ and $\theta > \theta_{1}$. This completes the proof.\\

\noindent\textbf{Remark 2.1} It is to be noted that stability of the predator-free equilibrium does not depend on the fractional-order $m$.\\

\noindent\textbf{Theorem 2.5} \textit{The following statements are true for the stability of the interior equilibrium point $E^{*}$ of system (\ref{Ecological fractional order model}).}
\textit{\begin{itemize}
		\item[(a)] If ~$trace(J^*) < 0$, i.e. if ~~$c_{2} < c < c_{1}$ with $\theta > \theta_{2}$, $\alpha> \frac{1}{Kh}$ then the interior equilibrium $E^{*}$ is locally asymptotically stable for $0< m \leq 1$, where $c_{2} = 1 - \frac{\theta + hd}{\alpha Kh(\theta - hd)}$, $c_{1} = 1 - \frac{d}{\alpha K (\theta - h d)}$ and $\theta_{2} = \frac{hd(\alpha Kh+1)}{\alpha Kh - 1}.$
		\item[(b)] If ~$0< trace(J^*) < 2\sqrt{det(J^*)}$, i.e. if ~~$0 < c < c_{2}$ with $\theta > \theta_{2}$, $\alpha> \frac{1}{Kh}$ then for any $m\in (0,m^*)$, the interior equilibrium $E^{*}$ is locally asymptotically stable and unstable for any $m\in (m^*,1]$. A Hopf bifurcation occurs at $m=m^*$, where $m^* = \frac{2}{\pi}\mid \cos^{-1}(\frac{trace(J^*)}{2\sqrt{det(J^*)}}) \mid$.
		\item[(c)] If ~$trace(J^*) \geq 2\sqrt{det(J^*)}$, then the interior equilibrium $E^{*}$ is unstable for any ~$m\in(0,1]$.
\end{itemize}}

\noindent\textbf{Proof}  For the interior equilibrium $E^{*}$, the Jacobian matrix is given by
\[\mathbf{J(E^{*})} = \begin{pmatrix}
r(1-\frac{2x^{*}}{ K}) - \frac{\alpha(1-c)y^{*}}{(1+\alpha(1-c)hx^{*})^{2}} & -\frac{\alpha(1-c)x^{*}}{(1+\alpha(1-c)hx^{*})} \\ \frac{\theta \alpha(1-c)y^{*}}{(1+\alpha(1-c)hx^{*})^{2}} & 0
\end{pmatrix}.\]
The corresponding characteristic equation is given by
\begin{equation}\label{chareqn}
\xi^{2} - trace(J^{*})\xi + det(J^{*})=0,
\end{equation}
where $trace(J^{*}) = \frac{rd\{\alpha hk(1-c)-1-2\alpha h(1-c)x^{*}\}}{\alpha \theta K(1-c)}$ and $det(J^{*}) = \frac{rd\{\alpha K(1-c)(\theta - hd)-d\}}{K\{\alpha(\theta - hd) + \alpha hd\}(1-c)}$. \\
Therefore, the roots of this equation are given by $$\xi_{1, 2} = \frac{1}{2}[trace(J^{*}) \pm \sqrt{trace(J^{*})^{2} - 4 det(J^{*})}].$$
\textit{(a)} Note that $trace(J^*)$ will be negative if $\alpha hk(1-c)-1-2\alpha h(1-c)x^{*} < 0$, i.e., if $c > 1 - \frac{\theta + hd}{\alpha Kh(\theta - hd)} = c_{2}$ with $\theta > \frac{hd(1+\alpha Kh)}{\alpha Kh - 1}$, $\alpha> \frac{1}{Kh}$. Since $det(J^{*}) > 0$, both roots of (\ref{chareqn}) are negative real or complex conjugate with negative real parts. Hence $\mid arg(\xi_{1, 2})\mid >\frac{m\pi}{2}, \forall m\in (0,1]$. So the positive interior equilibrium $E^{*}$ is locally asymptotically stable for $0<m \leq 1$ if $c_{2} < c < c_{1}$ with $\theta > max[hd + \frac{d}{\alpha K}, \frac{hd(1+\alpha Kh)}{\alpha Kh - 1}]$, $\alpha> \frac{1}{Kh}$. Noting that $max[hd + \frac{d}{\alpha K}, \frac{hd(1+\alpha Kh)}{\alpha Kh - 1}] = \frac{hd(1+\alpha Kh)}{\alpha Kh - 1} = \theta_{2}$, one gets the required stability result. This completes the proof of (a).\\

\textit{(b)} The condition $0< trace(J^*) < 2\sqrt{det(J^*)}$ will hold if $0 < c < c_{2}$ with $\theta > \theta_{2}$, $\alpha>\frac{1}{Kh}$, where $c_{2} = 1 - \frac{\theta + hd}{\alpha Kh(\theta - hd)}$. Since $0< trace(J^*) < 2\sqrt{det(J^*)}$, the equation (\ref{chareqn}) has two complex conjugate roots with positive real part given by
\begin{equation}\label{eigenvalue}
\xi_{i} = \frac{1}{2}[trace(J^{*}) \pm i\sqrt{4 det(J^{*}) - trace(J^{*})^{2}}], ~i=1,2,
\end{equation}
with $\mid arg(\xi_{1, 2})\mid  = \mid \cos^{-1}(\frac{trace(J^*)}{2\sqrt{det(J^*)}})\mid$. Assume that there exists a $m^* \in (0,1]$ such that $\mid \cos^{-1}(\frac{trace(J^*)}{2\sqrt{det(J^*)}}) \mid  = \frac{m^* \pi}{2}$. Then, following Theorem $3.1$, we have $\mid arg(\xi_{i})\mid >\frac{m\pi}{2}$ for all $m \in (0, m^*)$ and $\mid arg(\xi_{i})\mid <\frac{m\pi}{2}$ for all $m\in (m^*,1]$. Therefore, the positive interior equilibrium $E^{*}$ is locally asymptotically stable for ~$0< m < m^*$ and unstable for $m^* < m \leq 1$ when $0< trace(J^*) < 2\sqrt{det(J^*)}$, i.e., when $0 < c < c_{2}$.
A Hopf bifurcation will occur at $m=m^*$ under the following conditions \cite{Abdelouahab11,Xiang14}:
\begin{eqnarray} \nonumber
&(i)& Real(\xi_{i}) > 0,\nonumber\\
&(ii)& min_{i}\mid arg(\xi_{i})\mid = \frac{m^*\pi}{2}, ~i=1,2 \nonumber \\
&(iii)& \frac{d}{dm}[Real{\xi_i}]\mid_{m = m^*} \neq 0~ (transversality~~condition)\nonumber.
\end{eqnarray}
Note that $Real(\xi_{i})=\frac{1}{2}trace(J^*)>0$ and $\mid arg(\xi_{i})\mid = \frac{m^*\pi}{2}, ~i=1,2,$ by assumption. Also, $\frac{d}{dm}[Real{\xi_i}]\mid_{m = m^*}=\frac{\pi}{2}\neq 0.$ Therefore, a Hopf bifurcation exists as $m$ crosses the critical value $m^*$. The equilibrium $E^*$ is thus stable for all $m\in (0,m^*)$ and unstable for all $m\in (m^*,1].$ This completes the proof. \\

\textit{(c)} As $trace(J^*) \geq 2\sqrt{det(J^*)}$, the equation (\ref{chareqn}) has two real roots given by $\xi_{1, 2} = \frac{1}{2}[trace(J^{*}) \pm \sqrt{trace(J^{*})^{2} - 4 det(J^{*})}]$. Now for the positive root $\xi_{1}$, we note that $\mid arg(\xi_{1})\mid  = 0$. Since the eigenvalue $\xi_{1}$ does not satisfy $\mid arg(\xi_{1})\mid >\frac{m\pi}{2}, \forall m \in (0,1]$, therefore $E^{*}$ is  unstable for any $m \in (0,1]$. This completes the proof of (c). \\

\noindent\textbf{Theorem 2.6} \textit{The interior equilibrium $E^{*}$ is globally asymptotically stable for any $m\in (0, 1]$ if $c_{2} < c < c_{1}$ with $\theta > \theta_{2}, \alpha> \frac{1}{Kh}$, where $c_{2} = 1 - \frac{\theta + hd}{\alpha Kh(\theta - hd)}, ~ c_{1} = 1 - \frac{d}{\alpha K (\theta - h d)}~and~\theta_2 = \frac{hd(1+\alpha Kh)}{\alpha Kh - 1}$.} \\

\noindent\textbf{Proof}  Let us consider the Lyapunov function
$$V(x,y) = \bigg(x - x^{*} - x^{*}ln\frac{x}{x^{*}}\bigg) + \frac{1}{\theta -hd}\bigg(y - y^{*} - y^{*}ln\frac{y}{y^{*}}\bigg).$$
Here $V(x,y) > 0$ for all values of $x, y > 0$ and $V = 0$ only at $E^{*} = (x^{*},y^{*})$. Considering the $mth$ order fractional derivative of $V(x,y)$ along the solution of (\ref{Ecological fractional order model}) and using Lemma $2.4$, we have
\begin{eqnarray}\label{Globality_2}\nonumber
\begin{split}
^c_{0}D^{m}_{t}V(x,y) \leq & \frac{(x - x^{*})}{x} {^c_{0}}D^{m}_{t}x(t)+\frac{1}{\theta -hd}\frac{(y - y^{*})}{y} {^c_{0}}D^{m}_{t}y(t)\\
=& (x - x^{*})(r[1-\frac{x}{K}] - \frac{\alpha (1-c) y}{1+\alpha (1-c)h x}) \\
&+ \frac{\theta}{\theta-hd}(y - y^{*})(\frac{\alpha (1-c) x}{1+\alpha (1-c)h x} - \frac{d}{\theta})\\
=& (x - x^{*})[\frac{r(x^{*} - x)}{K} + \frac{\alpha (1-c) y^{*}}{1+\alpha (1-c)h x^{*}} - \frac{\alpha (1-c) y}{1+\alpha (1-c)h x}] \\
& + \frac{\theta}{\theta-hd}(y - y^{*})[\frac{\alpha (1-c) x}{1+\alpha (1-c)h x} - \frac{\alpha (1-c) x^{*}}{1+\alpha (1-c)h x^{*}}] \\
=& -\frac{r}{K}(x - x^{*})^{2}  + (x - x^{*}) \alpha (1-c)[\frac{y^{*}}{1+\alpha (1-c)h x^{*}} - \frac{y}{1+\alpha (1-c)h x}]\\
& + \frac{\theta}{\theta-hd}(y - y^{*}) \alpha (1-c)[\frac{x}{1+\alpha (1-c)h x} - \frac{x^{*}}{1+\alpha (1-c)h x^{*}}]\\
=& -\frac{r}{K}(x - x^{*})^{2}  - \frac{\alpha (1-c)(x - x^{*})(y - y^*)(1+\alpha (1-c)h x^{*})}{(1+\alpha (1-c)h x^{*})(1+\alpha (1-c)h x)} \\
& + \frac{\alpha^{2} (1-c)^{2} h y^{*}(x - x^{*})^2}{(1+\alpha (1-c)h x^{*})(1+\alpha (1-c)h x)} + \alpha (1-c)(x - x^{*})(y - y^*)\\
& \frac{(1+\alpha (1-c)h x^{*})}{(1+\alpha (1-c)h x^{*})(1+\alpha (1-c)h x)}\\
\leq & -\frac{r}{K}(x - x^{*})^{2}  +  \frac{\alpha^{2} (1-c)^{2} hy^{*}}{(1+\alpha (1-c)h x^{*})}(x - x^{*})^2\\
= & [\frac{r \alpha  h(1-c)(K - x^{*})}{K} - \frac{r}{K}](x - x^{*})^{2} \\
= & \frac{r}{K}[\alpha  h(1-c)(K-x^*) - 1] (x - x^{*})^{2}.
\end{split}
\end{eqnarray}
One can note that $^c_{0}D^{m}_{t}V(x,y) \leq0, \forall (x,y)\in R^{2}_{+}$ if $\alpha  h(1-c)(K-x^*) - 1 < 0$, i.e., if  $\alpha hK(1-c) < \frac{\theta}{\theta - hd} < \frac{\theta + hd}{\theta - hd}$. This implies $^c_{0}D^{m}_{t}V(x,y) \leq0, \forall (x,y)\in R^{2}_{+}$ if $c_{2} < c < c_{1}$, $\theta > \theta_{2}$, $\alpha> \frac{1}{Kh}$ and $^c_{0}D^{m}_{t}V(x,y) = 0$ implies that $(x,y) = (x^{*}, y^{*})$. Therefore, the only invariant set on which $^c_{0}D^{m}_{t}V(x,y) = 0$  is the singleton $\{E^{*}\}$. Then, following Lemma $4.6$ in \cite{HET15}, the interior equilibrium $E^{*}$ is globally asymptotically stable if the conditions in the theorem are satisfied. This completes the proof.
\section{Discretized fractional-order model and its analysis}
We first construct the discrete fractional-order model corresponding to the system (\ref{Ecological fractional order model}). Following Elsadany and Matouk \cite{Matouk15}, discretization of the model system (\ref{Ecological fractional order model}) with piecewise constant arguments can be done in the following manner:
\begin{eqnarray}\nonumber
^{c}_{0} D^{m}_{t}x & = & rx([t/s]s)\bigg(1-\frac{x([t/s]s)}{K}\bigg) - \frac{\alpha (1-c) x([t/s]s) y([t/s]s)}{1+\alpha (1-c)h x([t/s]s)}, \\
^{c}_{0} D^{m}_{t}y & =  & \frac{\theta \alpha (1-c) x([t/s]s) y([t/s]s)}{1+\alpha (1-c)h x([t/s]s)} - dy([t/s]s),\nonumber
\end{eqnarray}
with initial condition $x(0) = x_0 > 0$ and $y(0) = y_0> 0$.\\
Let $t\in [0,s)$, so that $t/s \in [0,1)$. In this case, we have
\begin{eqnarray}\nonumber
^{c}_{0} D^{m}_{t}x & = & x_0\bigg(r(1-\frac{x_0}{K}) - \frac{\alpha (1-c) y_0}{1+\alpha (1-c)h x_0}\bigg), \nonumber\\
^{c}_{0} D^{m}_{t}y & =  & y_0\bigg(\frac{\theta \alpha (1-c) x_0}{1+\alpha (1-c)h x_0} - d\bigg),\nonumber
\end{eqnarray}
and the solution of this fractional differential equation can be written as
\begin{eqnarray}\nonumber
x_1(t) & = & x_0 + J^{m}_{0}\bigg(x_0\bigg(r(1-\frac{x_0}{K}) - \frac{\alpha (1-c) y_0}{1+\alpha (1-c)h x_0}\bigg)\bigg) \nonumber \\
&  = & x_0 + \frac{t^m}{m\Gamma (m)} \bigg(x_0\bigg(r(1-\frac{x_0}{K}) - \frac{\alpha (1-c) y_0}{1+\alpha (1-c)h x_0}\bigg)\bigg), \nonumber \\
y_1(t) & = & y_0 + J^{m}_{0}\bigg(y_0\bigg(\frac{\theta \alpha (1-c) x_0}{1+\alpha (1-c)h x_0} - d\bigg)\bigg) \nonumber \\
& = & y_0 + \frac{t^m}{m \Gamma (m)}\bigg(y_0\bigg(\frac{\theta \alpha (1-c) x_0}{1+\alpha (1-c)h x_0} - d\bigg)\bigg).\nonumber
\end{eqnarray}
In the second step, we assume $t\in [s,2s)$ so that $t/s \in [1,2)$ and obtain
\begin{eqnarray}\nonumber
^{c}_{0} D^{m}_{t}x & = & x_1(s)\bigg(r(1-\frac{x_1(s)}{K}) - \frac{\alpha (1-c) y_1(s)}{1+\alpha (1-c)h x_1(s)}\bigg), \nonumber\\
^{c}_{0} D^{m}_{t}y & =  & y_1(s)\bigg(\frac{\theta \alpha (1-c) x_1(s)}{1+\alpha (1-c)h x_1(s)} - d\bigg).\nonumber
\end{eqnarray}
The solution of this equation reads
\begin{eqnarray}\nonumber
x_2(t) & = & x_1(s) + J^{m}_{s}\bigg(x_1(s)\bigg(r(1-\frac{x_1(s)}{K}) - \frac{\alpha (1-c) y_1(s)}{1+\alpha (1-c)h x_1(s)}\bigg)\bigg) \nonumber\\
& = & x_1(s) + \frac{(t - s)^m}{m\Gamma (m)}\bigg(x_1(s)\bigg(r(1-\frac{x_1(s)}{K}) - \frac{\alpha (1-c) y_1(s)}{1+\alpha (1-c)h x_1(s)}\bigg)\bigg), \nonumber\\
y_2(t) & = & y_1(s) + J^{m}_{s}\bigg(y_1(s)\bigg(\frac{\theta \alpha (1-c) x_1(s)}{1+\alpha (1-c)h x_1(s)} - d\bigg)\bigg) \nonumber\\
& = & y_1(s) + \frac{(t - s)^m}{m\Gamma (m)}\bigg(y_1(s)\bigg(\frac{\theta \alpha (1-c) x_1(s)}{1+\alpha (1-c)h x_1(s)} - d\bigg)\bigg),\nonumber
\end{eqnarray}
where $J^m_{s} = \frac{1}{\Gamma(m)}\int^{t}_{s} (t- \tau)^{(m-1)} d\tau, m>0$.\\ Repeating the discretization process $n$ times, we have
\begin{eqnarray}\nonumber
x_{n+1}(t) & = & x_n(ns) + \frac{(t - ns)^m}{m\Gamma (m)}\bigg(x_n(ns)\bigg(r(1-\frac{x_n(ns)}{K}) - \frac{\alpha (1-c) y_n(ns)}{1+\alpha (1-c)h x_n(ns)}\bigg)\bigg), \nonumber \\
y_{n+1}(t) & = & y_n(ns) + \frac{(t - ns)^m}{m\Gamma (m)}\bigg(y_n(ns)\bigg(\frac{\theta \alpha (1-c) x_n(ns)}{1+\alpha (1-c)h x_n(ns)} - d\bigg)\bigg), \nonumber
\end{eqnarray}
where $t\in [ns, (n+1)s)$.\\ Making $t\rightarrow (n+1)s$, we obtain the corresponding fractional discrete model of the continuous fractional model (\ref{Ecological fractional order model}) as
\begin{eqnarray}\nonumber \label{Discrete model}
x_{n+1} & = & x_n + \frac{s^m}{m\Gamma (m)}\bigg(x_n\bigg(r(1-\frac{x_n}{K}) - \frac{\alpha (1-c) y_n}{1+\alpha (1-c)h x_n}\bigg)\bigg), \nonumber \\
y_{n+1} & = & y_n + \frac{s^m}{m\Gamma (m)}\bigg(y_n\bigg(\frac{\theta \alpha (1-c) x_n}{1+\alpha (1-c)h x_n} - d\bigg)\bigg).
\end{eqnarray}
It is noticeable that Euler discrete model is a special case of this generalized discrete model when $m \rightarrow 1$.
\subsection{\textbf{Existence and stability of fixed of points}}
In the following, we investigate the dynamics of the discretized fractional-order model (\ref{Discrete model}). At the fixed point, we have $x_{n+1} = x_{n} = x$ and $y_{n+1} = y_{n} = y$. One can easily compute that (\ref{Discrete model}) has the same fixed points as in the fractional-order system (\ref{Ecological fractional order model}) given by $E_{0} = (0,0)$, $E_{1} = (K,0)$ and $E^{*} = (x^{*},y^{*})$, where
\begin{eqnarray}\nonumber
x^{*} & = & \frac{d}{\alpha (1-c)(\theta - h d)}, ~y^{*} = \frac{r(K-x^*)\{1 + \alpha h(1-c)x^{*}\}}{\alpha K(1-c)}.
\end{eqnarray}
The fixed point $E^{*}$ exists if $0 < c < c_{1}$ and $\theta > \theta_{1}$, where $c_{1} = 1 - \frac{d}{\alpha K (\theta - h d)}$, $\theta_{1} = hd + \frac{d}{\alpha K}$.\\

The Jacobian matrix of system (\ref{Discrete model}) at any arbitrary fixed point point $(x,y)$ reads
\begin{equation}
\mathbf{J(x,y)}=
\begin{pmatrix}
a_{11} & a_{12} \\ a_{21} & a_{22}
\end{pmatrix}
\label{eq:myeqn}
\end{equation}
where
\begin{eqnarray}\nonumber
a_{11} & =& 1 + \frac{s^m}{m\Gamma (m)}\bigg(r(1-\frac{x}{K}) - \frac{\alpha (1-c) y}{1+\alpha (1-c)h x}\bigg)\nonumber\\
&  & + \frac{s^m}{m\Gamma (m)} x\bigg(-\frac{r}{K} + \frac{\alpha^{2} h (1-c)^2 y}{(1+\alpha (1-c)h x)^2}\bigg), \nonumber \\
a_{12} & = & - \frac{s^m}{m\Gamma (m)} \frac{\alpha (1-c) x}{1 + \alpha (1-c) h x}, \nonumber \\
a_{21} & = & \frac{s^m}{m\Gamma (m)} \frac{\theta \alpha (1-c) y}{(1 + \alpha h (1-c) x)^2},\nonumber \\
a_{22} & =& 1 + \frac{s^m}{m\Gamma (m)} \bigg(\frac{\theta \alpha (1-c) x}{1+\alpha (1-c)h x} - d\bigg).\nonumber
\end{eqnarray}
Let $\xi_1$ and $\xi_2$ be the eigenvalues of the Jacobian matrix (\ref{eq:myeqn}). Then we have the following definition and lemma.\\

\noindent\textbf{Definition 3.1} \cite{Elaydi08,Milan15} \textit{A fixed point $(x, y)$ of system (\ref{Discrete model}) is called stable if $\mid\xi_1\mid < 1$, $\mid\xi_2\mid < 1$ and a source if $\mid\xi_1\mid > 1$, $\mid\xi_2\mid > 1$. It is called a saddle if $\mid\xi_1\mid < 1$, $\mid\xi_2\mid > 1$ or $\mid\xi_1\mid > 1$, $\mid\xi_2\mid < 1$ and a nonhyperbolic fixed point if either $\mid\xi_1\mid = 1$ or $\mid\xi_2\mid = 1$. It is called a spiral source if $\xi_{1,2} = \alpha \pm i \beta, \beta \neq 0, \alpha, \beta \in R$ and $\mid\xi_{1,2}\mid > 1$.}\\

\noindent\textbf{Lemma 3.1} \cite{Elaydi08} \textit{Let $\xi_1$ and $\xi_2$ be the eigenvalues of Jacobian matrix (\ref{eq:myeqn}). Then $\mid\xi_1\mid < 1$ and $\mid\xi_2\mid < 1$ if the following condition holds: $$(i) 1-det(J) > 0, (ii) 1-trace(J)+ det(J) > 0,~ and ~(iii) 1+ trace(J)+det(J) > 0.$$}
\noindent\textbf{Theorem 3.1} \textit{ (a) The fixed point $E_0$ is always unstable for $0<m\leq 1$. It will be a saddle point if $0 < s < s_1$ and a source if $s > s_1$. If $s = s_1$, then $E_0$ is nonhyperbolic, where $s_1 = \sqrt[m]{\frac{2m\Gamma(m)}{d}} $. } \\

\noindent \textit{(b) The fixed point $E_1$ is stable for $0<m\leq 1$ if $c > c_1$, $s < min \{s_2 , s_3\}$, where $s_2 = \sqrt[m]{\frac{2m\Gamma(m)}{r}}$, $s_3 = \sqrt[m]{\frac{2m\Gamma(m)\{1 + \alpha K h(1-c)\}}{d - K\alpha(1-c)(\theta- hd)}}$. It is a saddle point if $c > c_1$, $s_3 < s < s_2$~; or $c > c_1$, $s_2 < s < s_3$; and a source if $c > c_1$, $s > max\{s_2, s_3\}$.} \\

\noindent \textit{(c) The fixed point $E^*$ is locally asymptotically stable for $0<m\leq 1$ if $c_2 < c < c_1$ with $\theta > \theta_{2}, \alpha> \frac{1}{Kh}$ and $s < min\{s_4, s_5 \}$, where  $s_4 = \sqrt[m]{\frac{m\Gamma(m)G}{H}}$, $s_5 = \sqrt[m]{\frac{2m\Gamma(m)}{G}}$, \\
	$G = \frac{rx^*}{K\theta}[\theta +hd - \alpha hK (1-c)(\theta - hd)]$ and $H = \frac{rx^* (\theta - hd)}{K \theta }[\alpha K(1-c)(\theta - hd) -d]$}.\\

\noindent\textbf{Proof}  At the fixed point $E_0$,
the eigenvalues are $\xi_1 = 1 + r \frac{s^m}{m\Gamma (m)}$ and $\xi_2 = 1 - d \frac{s^m}{m\Gamma (m)}$. Since $|\xi_1| > 1$, $E_0$ is always unstable for $0<m \leq 1$. In fact, it is a saddle point if $0 < s < \sqrt[m]{\frac{2m\Gamma(m)}{d}}$ for which $|\xi_2 | < 1$ and a source if $s > \sqrt[m]{\frac{2m\Gamma(m)}{d}}$ for which $|\xi_2| > 1$. Again, it becomes nonhyperbolic if $s = \sqrt[m]{\frac{2m\Gamma(m)}{d}}$ for any $m \in (0,1]$. \\

\noindent The eigenvalues evaluated at the fixed point $E_1$ are evaluated as

$$\xi_1 = 1 - r \frac{s^m}{m\Gamma (m)}, ~\xi_2 = 1 + \frac{s^m}{m\Gamma (m)} \bigg(\frac{\theta \alpha (1-c) K}{1+\alpha (1-c)h K} - d\bigg).$$
Note that for $0<m \leq 1$,  $|\xi_{1,2} | < 1$ hold if $$s < min \{\sqrt[m]{\frac{2m\Gamma(m)}{r}},  \sqrt[m]{\frac{2m\Gamma(m)\{1 + \alpha K h(1-c)\}}{d - K\alpha(1-c)(\theta- hd)}}\}.$$ Therefore, $E_1$ is locally asymptotically stable for $0<m\leq 1$ if $c > c_1$ and $s < min \{s_2, s_3\}$. However, $|\xi_1| > 1$ if $s > s_2$ and $|\xi_2| > 1$ if $s > s_3$ with $c>c_1$. Thus, $E_1$ will be a source if $c>c_1$ and $s > max\{s_2, s_3\}$. The fixed point $E_1$ will be a saddle point if either of the conditions $(i)~ s_3< s < s_2 $, $c > c_1$ or $(ii)~ s_2< s < s_3 $, $c > c_1$ holds.\\

\noindent At the interior fixed point $E^*$, the Jacobian matrix is evaluated as
\[
\mathbf{J(x^*,y^*)} = \begin{pmatrix}
a_{11} & a_{12} \\ a_{21} & a_{22} \end{pmatrix},\]
where $a_{11} = 1 - \frac{s^m}{m\Gamma (m)}G,$ $a_{12} = - \frac{s^m}{m\Gamma (m)} \frac{\alpha (1-c) (\theta -hd )x^*}{\theta},$
$a_{21}  =  \frac{s^m}{m\Gamma (m)} \frac{r (\theta - h d) (K - x^*)}{K},$ $a_{22}  = 1$
and $a_{12}a_{21} = - {(\frac{s^m}{m\Gamma (m)})}^{2} H$ with $G = \frac{rx^*}{K\theta}[\theta +hd - \alpha hK (1-c)(\theta - hd)]$ and \\
$H = \frac{rx^* (\theta - hd)}{K \theta }[\alpha K(1-c)(\theta - hd) -d]$.

\noindent Note that $H > 0$ if $c < c_1$ and $G > 0$ if $c>c_2, \theta > \theta_2, \alpha > \frac{1}{Kh}$. After some algebraic manipulations, we have $$det(J) = 1- (\frac{s^m}{m\Gamma (m)})G + {(\frac{s^m}{m\Gamma (m)})}^2 H ~\mbox{and}~ trace(J) = 2 - (\frac{s^m}{m\Gamma (m)}) G.$$
Thus, $1 - trace(J) + det(J) = {(\frac{s^m}{m\Gamma (m)})}^2 H > 0$ if $c< c_1$. Also, $1-det(J) = (\frac{s^m}{m\Gamma (m)}) (G - (\frac{s^m}{m\Gamma (m)}) H)$ is positive if $s < s_4$, where $s_4 = \sqrt[m]{\frac{m\Gamma(m)G}{H}}$ and $c> c_2$ with $\theta > \theta_{2}, \alpha> \frac{1}{Kh}$.

\noindent One can compute that
$$1 + trace(J) + det(J) = 2(2 - (\frac{s^m}{m\Gamma (m)})G) + {(\frac{s^m}{m\Gamma (m)})}^2 H.$$ \\
This expression will be positive if $0< s< s_5$, where $s_5 = \sqrt[m]{\frac{2m\Gamma(m)}{G}}$. Therefore, the fixed point $E^*$ is stable if $c_2 < c < c_1$ and $s < min\{s_4, s_5\}$ for any $m \in (0,1]$ and unstable otherwise. Hence the theorem.\\

\noindent\textbf{Remark 3.1} Here we also observe that the predator-free fixed point $E_1$ looses its stability through transcritical bifurcation (a real eigenvalue that passes through $+1$) when $1 - trace(J) + det(J) = 0$ at $c = c_1$ for any $m$ \cite{Matouk15}. Again our model system (\ref{Discrete model}) undergoes a flip bifurcation (a real eigenvalue becomes equal to $-1$) when $1+trace(J)+det(J) = 0$ at the predator-free fixed point $E_1$ for $c = c_1$ and $s = s_5 = \sqrt[m]{\frac{2 m\Gamma(m)}{G}}$.\\

\noindent\textbf{Remark 3.2} Note that the eigenvalues of $J(x^*,y^*)$ are $$\xi_{1,2} = \frac{1}{2}[2 - \frac{s^m}{m\Gamma(m)}G \pm \frac{s^m}{m\Gamma(m)} \sqrt{G^2 - 4H}].$$ Therefore, $\xi_{1,2}$ are complex conjugate if $G^2 -4H< 0$, i.e., if $-2\sqrt{H}< G < 2\sqrt{H}$.
\noindent Now, $$\mid \xi_{1,2}\mid = 1- (\frac{s^m}{m\Gamma (m)})G + {(\frac{s^m}{m\Gamma (m)})}^2 H = det(J)$$ and this modulus is equal to unity if $det(J) = 1$, i.e., if $s = \sqrt[m]{\frac{m\Gamma(m)G}{H}} = s_4$. Since $G = \frac{s^m}{m\Gamma(m)}H >0,$ the previous inequality becomes $0< G < 2\sqrt{H}$. Therefore, we can conclude that $J(x^*,y^*)$ has complex conjugate roots with unit modulus if parameters belong to the set $$U = \{(m, s, r,K, \alpha, \theta, h, c,d):  0< G < 2\sqrt{H}, s = \sqrt[m]{\frac{m\Gamma(m)G}{H}}\}.$$ Therefore, if the parameter $s$ varies in the neighborhood of $s_4$ and $(m, s, r,K, \alpha, \theta, h, c,d) \in U$, the system (\ref{Discrete model}) may undergo a Hopf bifurcation around the equilibrium $E^*$. \\

\subsection{Hopf Bifurcation and its stability}
Here we prove the existence of Hopf bifurcation around $E^{*} = (x^*, y^*)$ and its stability. Let $S_1 = \frac{s^m}{m\Gamma (m)}$ and $S^*$ be a perturbation in the bifurcation parameter $S_1$, where $|S^*|<<1$. Then a perturbation form of model (\ref{Discrete model}) can be represented as \cite{AET15}
\begin{eqnarray}\label{map_1}
x_{n+1} & = & x_n + (S_1 + S^*)\bigg(x_n\bigg(r(1-\frac{x_n}{K}) - \frac{\alpha (1-c) y_n}{1+\alpha (1-c)h x_n}\bigg)\bigg), \nonumber \\
y_{n+1} & = & y_n + (S_1 + S^*)\bigg(y_n\bigg(\frac{\theta \alpha (1-c) x_n}{1+\alpha (1-c)h x_n} - d\bigg)\bigg).
\end{eqnarray}
Let $X_n = x_n - x^*, Y_n = y_n -  y^*$ so that the fixed point $E^* = (x^*, y^*)$ of the map (\ref{map_1}) is transformed into the origin. The transformed system reads
\begin{eqnarray}\label{map_2}
X_{n+1} & = & c_{11}X_n + c_{12}Y_n + c_{13}X_n Y_n, \nonumber \\
Y_{n+1} & = & c_{21}X_n + c_{22}Y_n + c_{23}X_n Y_n,
\end{eqnarray}
where $c_{11} = 1 - (S_1 + S^*)G,$ $c_{12} = - (S_1 + S^*) \frac{\alpha (1-c) (\theta -hd )x^*}{\theta},$
$c_{21}  =  (S_1 + S^*) \frac{r (\theta - h d) (K - x^*)}{K},$ $c_{22}  = 1,$ $c_{13} = -\frac{\alpha (1-c)(S_1 + S^*)}{2(1 + \alpha(1-c)hx^*)^2},$ $c_{23} = \frac{\theta\alpha (1-c)(S_1 + S^*)}{2(1 + \alpha(1-c)hx^*)^2},$ $c_{12}c_{21} = -(S_1 + S^*)^2 H$ and $S_1 = \frac{s_4^m}{m\Gamma (m)}$ with $G = \frac{rx^*}{K\theta}[\theta +hd - \alpha hK (1-c)(\theta - hd)]$, $H = \frac{rx^* (\theta - hd)}{K \theta }[\alpha K(1-c)(\theta - hd) -d]$.\\

The characteristic equation associated with the linearization of the model (\ref{map_2})
at $(X_n, Y_n) = (0, 0)$ is given by
\begin{equation}\label{chareqn_1}
\lambda^2 + p(S^*)\lambda + q(S^*) = 0,
\end{equation}
where
\begin{eqnarray}\label{coeff_1}
p(S^*) = -2 + G(S_1 + S^*), ~~q(S^*) = 1 - G(S_1 + S^*) + H (S_1 + S^*)^2.
\end{eqnarray}
Since the parameters $(m, s, r,K, \alpha, \theta, h, c,d) \in U$  and $S^*$ varies in a small neighborhood of $S^* = 0$, and the roots of (\ref{chareqn_1}) are pair of complex conjugate numbers $\lambda_1$ and $\lambda_2$ denoted by
\begin{eqnarray}
\lambda_{1,2} &=& \frac{-p(S^*) \pm i \sqrt{4q(S^*) - p^2(S^*)}}{2},\nonumber \\
&=& \frac{1}{2}[2 - (S_1 + S^*)G \pm i(S_1 + S^*) \sqrt{4H - G^2}].
\end{eqnarray}
Therefore, $|\lambda_{1,2}| = \sqrt{q(S^*)}$. Since $q(S^*) = 1$ at $S^* = 0$, when $s = s_4 = \sqrt[m]{\frac{m\Gamma(m)G}{H}}$, then $|\lambda_{1,2}| = 1$ at $S^* = 0$ for $s = s_4$.\\
Consequently for $s = s_4$,
$$\frac{d|\lambda_{1,2}|}{dS^*}|_{S^* = 0} = \frac{G}{2} \neq 0 ~~\mbox{(transversality condition)}.$$
Also, at $S^* = 0$, $\lambda_{1,2}^n \neq 1$ for $n = 1,2,3,4$ (nonresonance conditions), which is equivalent to
$$p(0)~~\neq~~ -2,0,1,2.$$ Since $p^2(0) - 4q(0) < 0$ and $q(0) = 1$, we have $p^2(0) < 4$; then $p(0) \neq \pm2$. It is only require that $p(0) \neq 0, 1$, which leads to
\begin{equation} \label{cond_1}
G^2 \neq 3H,2H
\end{equation}
for $s = s_4.$
\noindent Next, we study the normal form of the model (\ref{map_2}) at $S^* = 0$. Let $\delta = Re(\lambda_{1,2})$ and $\beta = Im(\lambda_{1,2})$. We construct an invertible matrix
\begin{equation} \nonumber
{T}=
\begin{pmatrix}
c_{12} & 0 \\ \delta - c_{11} & -\beta
\end{pmatrix}
\end{equation}
and consider the translation
\begin{equation} \nonumber
\begin{pmatrix}
X_n  \\ Y_n
\end{pmatrix}
= T\begin{pmatrix}
u_n  \\ v_n
\end{pmatrix}.
\end{equation}
Thus, the map (\ref{map_2}) becomes
\begin{equation} \label{map_3}
\begin{pmatrix}
u_{n+1}  \\ v_{n+1}
\end{pmatrix}
\rightarrow \begin{pmatrix}
\delta & -\beta \\ \beta & \delta
\end{pmatrix}
\begin{pmatrix}
u_n  \\ v_n
\end{pmatrix} +
\begin{pmatrix}
P(u_n, v_n)  \\ Q(u_n, v_n)
\end{pmatrix},
\end{equation}
where
\begin{eqnarray}
P(u_n, v_n) &=& c_{13}[(\delta - c_{11})u_n^2 - \beta u_n v_n], \nonumber\\
Q(u_n, v_n)& =& ((c_{11} - \delta)c_{13} + c_{12}c_{23})[\frac{(c_{11}-\delta)}{\beta}u_n^2 + u_n v_n],\\
X_n &=& c_{12}u_n, \nonumber\\
Y_n &=& (\delta - c_{11})u_n - \beta v_n. \nonumber
\end{eqnarray}

\noindent In order to undergo Hopf bifurcation, we require that the following discriminatory quantity $\gamma$ be nonzero
\begin{equation}\label{cond_2}
\gamma = \bigg[-Re\bigg(\frac{(1-2\lambda_1)\lambda_2^2}{1-\lambda_1}\xi_{11}\xi_{20}\bigg)- \frac{1}{2}|\xi_{11}|^2 - |\xi_{02}|^2 + Re(\lambda_2\xi_{21})\bigg]\bigg|_{S^* = 0},
\end{equation}
where
$$\lambda_1 = \delta + i \beta, \lambda_2 = \delta - i \beta, \nonumber $$
$$\xi_{11} = \frac{1}{4}[(P_{u_nu_n} + P_{v_n v_n}) + i (Q_{u_nu_n} + Q_{v_n v_n})], \nonumber$$
$$\xi_{20} = \frac{1}{8}[(P_{u_nu_n} - P_{v_n v_n} + 2Q_{u_nv_n})+ i (Q_{u_nu_n} - Q_{v_n v_n} - 2P_{u_nv_n})], \nonumber $$
$$\xi_{02} = \frac{1}{8}([P_{u_nu_n} - P_{v_n v_n} - 2Q_{u_nv_n}) + i (Q_{u_nu_n} - Q_{v_n v_n} + 2P_{u_nv_n})],~~\nonumber $$
$$\xi_{21} = \frac{1}{16}[(P_{u_nu_nu_n} + P_{u_nv_n v_n} +Q_{u_nu_nv_n} + Q_{v_nv_nv_n}) + i (Q_{u_nu_nu_n} + Q_{u_nv_n v_n} \nonumber$$
$$- P_{u_nu_nv_n} - P_{v_nv_nv_n})], \nonumber $$
$$P_{u_nu_n} = 2c_{13}(\delta - c_{11}),~~ P_{v_n v_n} = 0,~~ P_{u_n v_n} = -\beta c_{13}, \nonumber$$
$$Q_{u_nu_n} = 2((c_{11} - \delta)c_{13} + c_{12}c_{23})\frac{(c_{11} - \delta)}{\beta},~~ Q_{v_n v_n} = 0,~~Q_{u_n v_n} = (c_{11} - \delta)c_{13} + c_{12}c_{23}, \nonumber$$
$$P_{u_nu_nu_n} = P_{u_nu_nv_n} = P_{u_nv_nv_n} = P_{v_nv_nv_n} = 0, \nonumber $$
$$Q_{u_nu_nu_n} = Q_{u_nu_nv_n} = Q_{u_nv_nv_n} = Q_{v_nv_nv_n} = 0, \nonumber $$

From the above analysis and Theorem $3.2$ in \cite{He11}, following theorem can be stated.\\

\noindent \textbf{Theorem 3.2} \textit{If conditions (\ref{cond_1}) and (\ref{cond_2}) hold, then the system (\ref{Discrete model}) undergoes Hopf bifurcation at the positive fixed point $E^{*}=(x^*, y^*)$ when the parameter $s$ varies in the small neighborhood of $s_4$. Furthermore if $\gamma < 0$ (respectively $\gamma > 0$), then an attracting (respectively repelling) invariant closed curve bifurcates from the fixed point $E^{*}=(x^*, y^*)$ for $s> s_4$ (respectively $s< s_4$), where $s_4 = \sqrt[m]{\frac{m\Gamma(m)G}{H}}$.}\\

A comparison table on dynamical behaviors of system (\ref{Ecological_model}) with corresponding fractional-order and discretized fractional-order versions has been given in Table $3.1$.\\

\noindent{\bf Table 3.1. Comparison of dynamical behaviors of three systems.}
\begin{flushleft}
	{\small
		\begin{tabular}
			{|l l l l|} \hline
			Equilibrium  & Continuous  & Fractional order & Discretized fractional \\
			point	& system \cite{Bairagi11} & system & order system \\  \hline
			
			$~~~~~E_0$    & Unstable    & Unstable  & Saddle if $0 < s < s_1$,   \\
			&                        &  for all $m \in (0, 1]$             & Source if $s > s_1$ $\&$
			\\
			&                        &              & Nonhyperbolic if $s = s_1$  \\
			&                        &              & for all $m \in (0, 1]$ \\
			\hline
			$~~~~~E_1$  & LAS if $c>c_{1}$, $\theta > \theta_{1}$  & Same with the  & LAS if $c > c_1$ and  \\
			&   $\&$ Unstable if $c<c_{1}$                     & continuous system  & $s < min \{s_2 , s_3\}$,    \\
			&                        & for all $m \in (0, 1]$  & Saddle if $c > c_1$, $s_3 < s < s_2$;  \\
			&                        &              & or $c > c_1$, $s_2 < s < s_3$,   \\
			&                        &              & Source if $c > c_1, s > max\{s_2, s_3\}$  \\
			&                         &             & for all $m \in (0, 1]$\\
			\hline			
			$~~~~~E^*$&LAS if $c_{2} < c < c_{1}$  & $(i)$ LAS if $c_{2} < c < c_{1}$ & LAS if $c_2 < c < c_1$ with $\theta > \theta_{2}, $ \\
			& with $\theta > \theta_{2}$, $\alpha> \frac{1}{Kh}$ $\&$ &  with $\theta > \theta_{2}$, $\alpha> \frac{1}{Kh}$              &  $\alpha> \frac{1}{Kh}$  and $s < min\{s_4, s_5 \}$    \\
			&  Unstable if $0 < c < c_{2}$   &   for all $m \in (0, 1]$               &  for all $m \in (0, 1]$   \\
			& with $\theta > max\{\theta_1, \theta_{2}\}$,  & $(ii)$ LAS if $0 < c < c_{2}$   &  \\
			& $\alpha> \frac{1}{Kh}$  &with $\theta > \theta_{2}$, $\alpha> \frac{1}{Kh}$             &  \\
			&   & for all $m \in (0, m^*)$                &  \\
			&   & $\&$ Unstable                           &  \\	
			&   &  for all $m \in (m^*, 1]$               &  \\					
			\hline
	\end{tabular}}
\end{flushleft}

\section{Numerical Simulations}
In this section, we perform extensive numerical computations of fractional-order differential equations (FDE) system (\ref{Ecological fractional order model}) for different fractional values of $m$ $(0 < m \leq 1)$ as well as the fractional-order discrete system (\ref{Discrete model}). We use Adams-type predictor corrector method for the numerical solution of FDE system (\ref{Ecological fractional order model}). It is an effective method to give numerical solutions of both linear and nonlinear FDE \cite{Diethelm02,Diethelm04}. We first replace our system (\ref{Ecological fractional order model}) by the following equivalent fractional integral equations:
\begin{eqnarray}\label{Eco-epidemiological fractional integral eqn} \nonumber
x(t) & = & x(0) + D^{-m}_{t} [rx\bigg(1-\frac{x}{K}\bigg) - \frac{\alpha (1-c) x y}{1+\alpha (1-c)h x}], \nonumber\\
y(t) & = & y(0) + D^{-m}_{t} [\frac{\theta \alpha (1-c) x y}{1+\alpha (1-c)h x} - dy],
\end{eqnarray}
and then apply the PECE (Predict, Evaluate, Correct, Evaluate) method. 

Three examples are presented to illustrate the analytical results of FDE system obtained in the previous section. To explore the effect of habitat complexity and fractional-order, we varied $c$ and $m$ in their respective ranges $0<c<1$ and $0< m<1$. We also plotted the solutions for $m=1$, whenever necessary, to compare the solutions of fractional-order system with that of integer order system.\\

\textbf{Example 1:} We considered the parameter values as $r = 2.65$, $K = 898$, $\alpha = 0.045$, $h = 0.0437$, $d = 1.06$ and initial point $x(0) = 10, y(0) = 5$ from \cite{Bairagi11}. Step size in all simulations is considered as $0.05$. Note that the condition $\alpha>\frac{1}{Kh}$ is always satisfied by this parameter set. Following Theorem $2.3(b)$, we compute  $c_{1} = 0.8445$, $\theta_{1} = 0.0726$ and select $c = 0.86 (> c_{1})$, $\theta = 0.215 (> \theta_{1})$ to show that the predator-free equilibrium $E_{1}$ of the system (\ref{Ecological fractional order model}) is asymptotically stable for all $m \in (0, 1]$ (Fig. 1). It is noticeable that the solutions reach to the equilibrium more slowly as the value of $m$ gets smaller. The phase planes presented in Fig. 2 show that the solution trajectory with different initial conditions (denoted by stars) reach to the equilibrium point (red circle) in each case, following Theorem $2.4$, depicting the global stability of the predator-free equilibrium $E_{1}$ for different values of $m$.

\begin{figure} [H] \label{time series_E_1}
	\hspace{-1in}	
	\includegraphics[width=6.5in, height=2in]{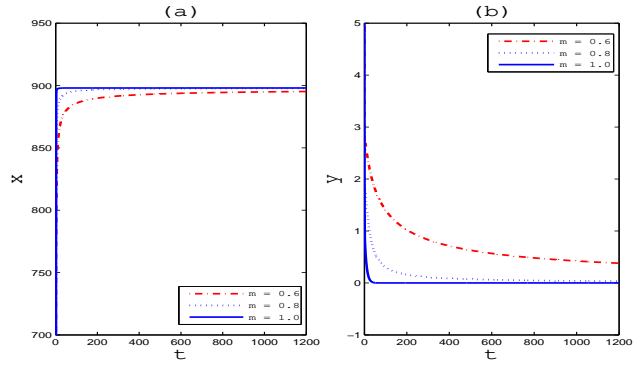}
	\caption{Asymptotically stable solutions of $x$ (prey) and $y$ (predator) for different fractional orders ($0<m<1$) and integer order $m = 1$ (solid line). It shows that the convergence rate of solutions to the equilibrium value is slower as $m$ becomes smaller. Here parameters are $r = 2.65$, $K = 898$, $\alpha = 0.045$, $h = 0.0437$, $d = 1.06$ and $c = 0.86,$ $\theta = 0.215$.}
	\label{Stable_1.eps}
\end{figure}

\begin{figure}[H] \label{global stability_E_1}
	\hspace{-1in}
	\includegraphics[width=6.5in, height=2in]{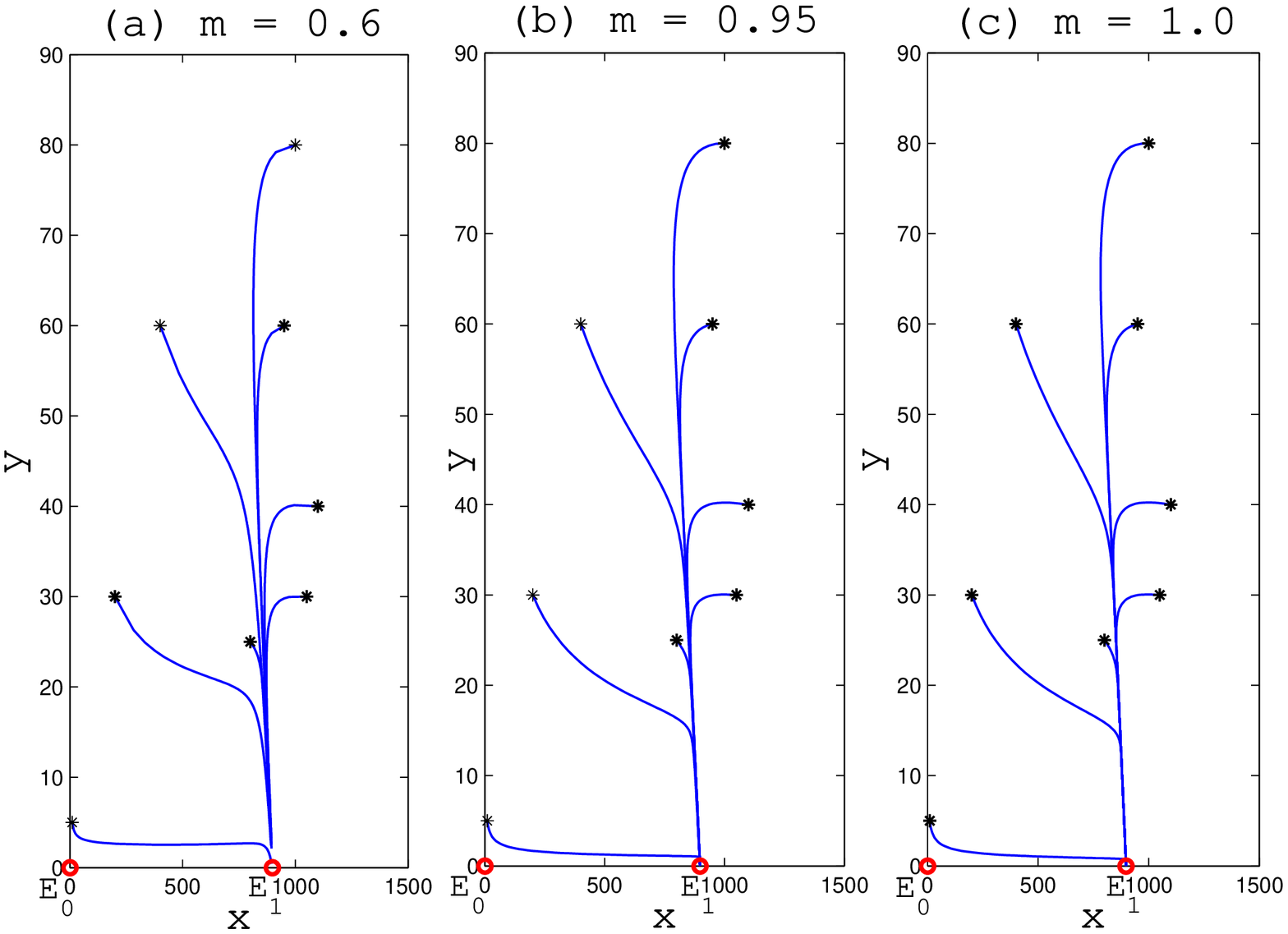}
	\caption{Trajectories with different initial values converge to the predator-free equilibrium $E_{1}$ for different values of $m$, indicating global stability of the equilibrium $E_{1}$, when conditions of Theorem 2.4 are satisfied. All parameters are as in Fig. 1.}
	\label{Stable_2.eps}
\end{figure}

\noindent \textbf{Example 2:} For the same parameter values as in Example 1, we compute $\theta_{2} = 0.1673$, $c_{2} = 0.1227$, $c_{1} = 0.8445$ and $trace(J(E^*)) = -0.3398 < 0$. Thus, following Theorem $2.5 (a)$, if we choose $\theta = 0.215 ~(> \theta_{2})$ and $c=0.45 ~(c_2<c<c_1)$ then solutions for all $m$ eventually converge to the equilibrium point $E^{*}$ where both the prey and predator populations coexist in the form of a stable equilibrium (Fig. 3).
\begin{figure}[H]
	\hspace{-1in}
	\includegraphics[width=6.5in, height=2.5in]{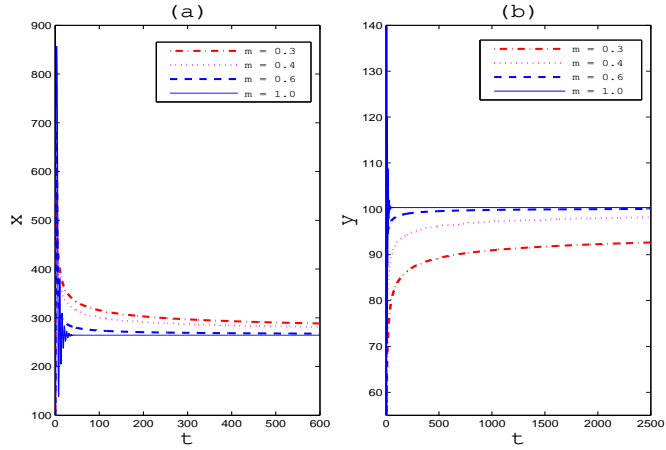}
	\caption{Asymptotically stable solutions of $x$ and $y$ populations for different fractional orders $m$ $(0<m<1)$ and integer order $m = 1$ (solid line). Parameters are as in Fig. 1 except $c = 0.45$.}
	\label{global_1.eps}
\end{figure}
If we choose $c<c_2=0.1227,$ say $c=0.05,$ and $\theta=0.215 (<\theta_2)$ as before then we obtain $trace(J^*) = 0.0437$, $2 \sqrt{det(J^*)} = 2.7152$. Therefore, from Theorem $2.5 (b)$, there exists a critical value $m=m^*=0.9898$ below which $E^*$ is stable and above which it is unstable. The stable behavior of the system (\ref{Ecological fractional order model}) for $m=0.95 (<0.9898)$ is presented in Fig. 4a and the unstable behavior of the system for $m=0.995 (>0.9898)$ in Fig. 4b. A Hopf bifurcation occurs at $m=m^*$. One can obtain a critical value $m^*$ for each $c\in (0, c_2)$, following Theorem 2.5(b), and can draw a stability region of $E^*$ in $c-m$ plane. The bifurcation curve separates the stable and unstable region (see Fig. 4c).
\begin{figure}[H]
	\includegraphics[width=8.5in, height=2.5in]{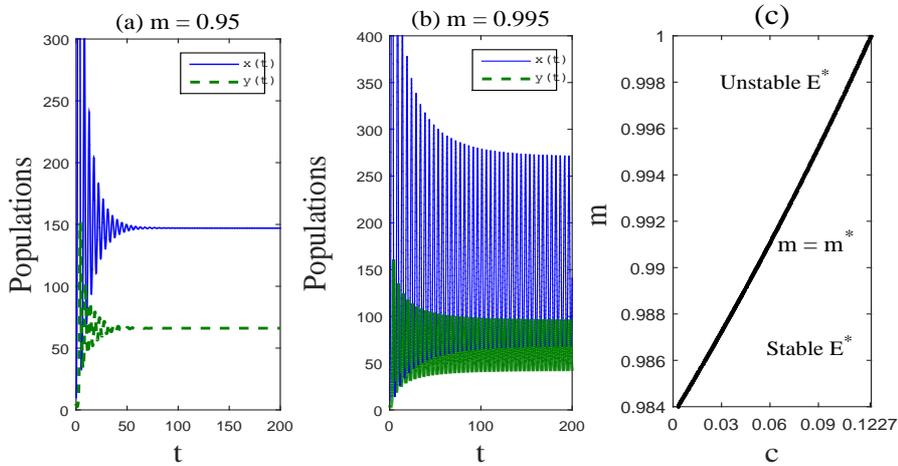}
	\caption{(a) Stable behavior of the system (\ref{Ecological fractional order model}) for $m=0.95$ ($<m^*=0.9898$), (b) unstable behavior of the system (\ref{Ecological fractional order model}) for $m = 0.995 (>m^*=0.9898)$. Here $c = 0.05$ and other parameters are as in Fig. 1. (c) Stability region of $E^*$ in $c-m$ plane when $c \in (0, c_2)$.}
	\label{stable_phase.eps}
\end{figure}
\noindent \textbf{Example 3:} Global stable behavior of system (\ref{Ecological fractional order model}) around the interior equilibrium $E^*$ is presented in Fig. 5. This figure shows that solutions with different initial conditions converge to the coexisting equilibrium $E^*$ for all values of $m$ when the conditions of Theorem $2.6$ are fulfilled.
\begin{figure}[H]
	\hspace{-1in}
	\includegraphics[width=6.5in, height=2.5in]{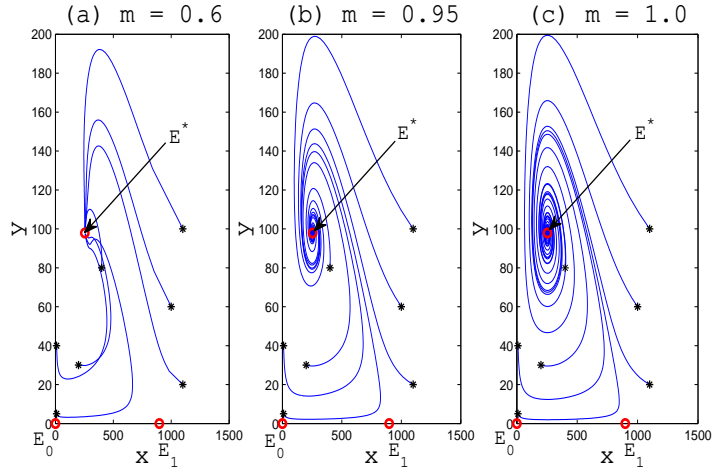}
	\caption{Global stability of the interior equilibrium $E^*$ for different values of $m$. Trajectories with different initial values converge to $E^{*}$ when the degree of habitat complexity is intermediate ($c_2 < c< c_1$), where $c_{2} = 0.1227$ and $c_{1} = 0.8445$. Here $c=0.45$ and other parameters are as in Fig. 1.}
	\label{limit_cycle_1.eps}
\end{figure}

\noindent \textbf{Example 4:} To illustrate the corresponding discrete system (\ref{Discrete model}) of the fractional-order system (\ref{Ecological fractional order model}), we consider the same parameter set as in Example 1. Stability of fixed points depends on the step size $s$ for different fractional-order $m$ (see Theorem 3.1). Assigning $c = 0.86$ for $E_1$ and $c = 0.45$ for $E^*$, we present the ranges of step size for the stability of the corresponding fixed point for different values of fractional-order $m$ in Table 5.1.\\

\noindent{\bf Table 5.1.} Restriction on the step size, following Theorem 3.1, for the stability of fixed points $E_1$ and $E^*$ for different fractional-order $m$.

\hspace{0.1in}
{\small
	\begin{tabular}
		{|l l l|} \hline
		& $E_1$ & $E^*$ \\ \hline
		Fractional order $m$ & Step size $s<min(s_2,s_3)$ & Step size $s<min(s_4,s_5)$  \\  \hline
		
		$m = 0.3$                 &  $s_2 = 0.2729$,         & $s_4 = 0.0041$, \\
		&  $s_3 = 26269$           & $s_5 = 256.7923$  \\
		\hline
		$m = 0.4$                 & $s_2 = 0.3669$,          & $s_4 = 0.0159$, \\
		&  $s_3 = 2005.2$           & $s_5 = 62.3401$   \\
		\hline
		$m = 0.6$                 & $s_2 = 0.5186$,           & $s_4 = 0.0639$, \\
		&  $s_3 = 160.8894$         & $s_5 = 15.9072$   \\
		\hline
		$m = 0.8$                 & $s_2 = 0.6436$,            & $s_4 = 0.1339$, \\
		&   $s_3 = 47.5805$          & $s_5 = 8.3894$   \\
		\hline
		$m = 0.95$                & $s_2 = 0.7279$,              & $s_4 = 0.1940$, \\
		&  $s_3 = 27.2757$            & $s_5 = 6.3253$   \\
		\hline
\end{tabular}}

\vspace{2cm}
\begin{figure}[H]
	\centering
	\includegraphics[width=10in, height=4in]{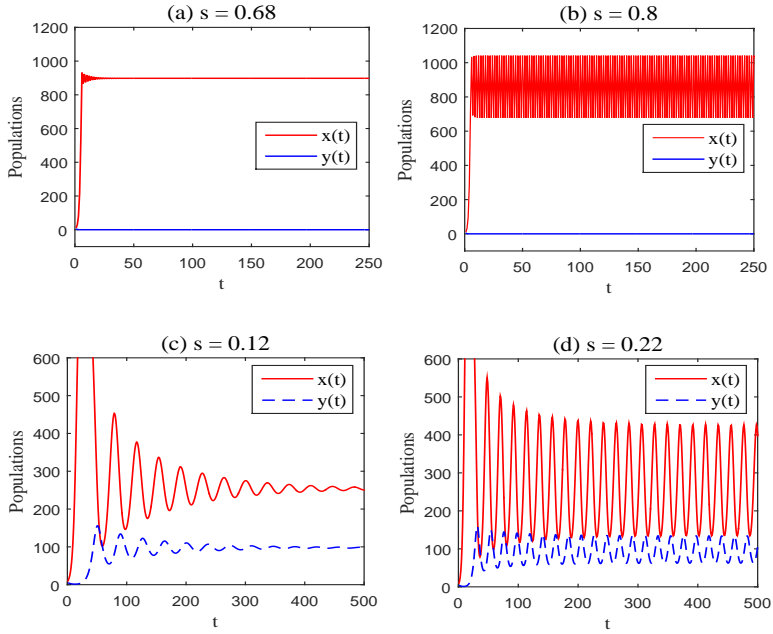}
	\vspace{-2cm}
	\caption{Stable and unstable behavior of different fixed points for different step sizes. Here step size has been considered as $s = 0.68$ and $s = 0.8$ to show the stability (Fig. a) and instability (Fig. b) of the fixed point $E_1$ when $c = 0.86$ and $m = 0.95$. Similar behavior for the fixed point $E^*$ has been shown for $s = 0.12$ (Fig. c) and $s = 0.22$ (Fig. d) when $c = 0.45$ and $m = 0.95$. Other parameters are as in Fig. 1.}
	\label{limit_cycle_2.eps}
\end{figure}
For example, when $m = 0.95$  then the step size $s$ should be less than $0.7279$ and $0.1940$ for $E_1$ and $E^*$, respectively, to be stable whenever they exist and unstable if it exceeds. In Fig. 6a, we have plotted the stable behavior of the fixed point $E_1$ for $s = 0.68$ and the unstable oscillatory behavior is presented in Fig. 6b for $s = 0.8$. Similar stable and unstable behaviors of the fixed point $E^*$ are presented in Figs. 6c-6d for $s = 0.12$ and $s = 0.22$, respectively. Following Theorem 3.2, we obtain a pair of complex conjugate eigenvalues as $\lambda_{1,2} =  0.9635 \pm 0.2678i$, where $\mid \lambda_{1,2}\mid = 1,~\frac{d|\lambda_{1,2}|}{dS^*}\bigg|_{S^* = 0} = 0.1699 > 0$ and $\gamma = -0.000 000019961 < 0$ at $(m,s,r,K,\alpha, \theta,h,c,d) = (0.95, 0.194, 2.65, 898, 0.045, 0.215, 0.0437, 0.45, 1.06) \in U$.  This implies that the system (\ref{Discrete model}) undergoes a Hopf bifurcation at the fixed point $E^* = ( 253.9056, 97.8867)$ for $s = 0.194~(=s_4)$. The bifurcation diagrams (Figs. 7a, 7c) represents it succiently. Also, the bifurcating closed curve is stable as $\gamma$ is negative. The system shows period doubling bifurcations leading to chaos as the step-size is further increased. Clear maginifed pictures of period doubling are presented in Figs. 7b \& 7d. Phase diagrams of the fractional-order discrete system (\ref{Discrete model}) for some particular values of $s$ are presented in Fig. 8.

\begin{figure}[H]
	\centering	
	\includegraphics[width=8.5in, height=3.5in]{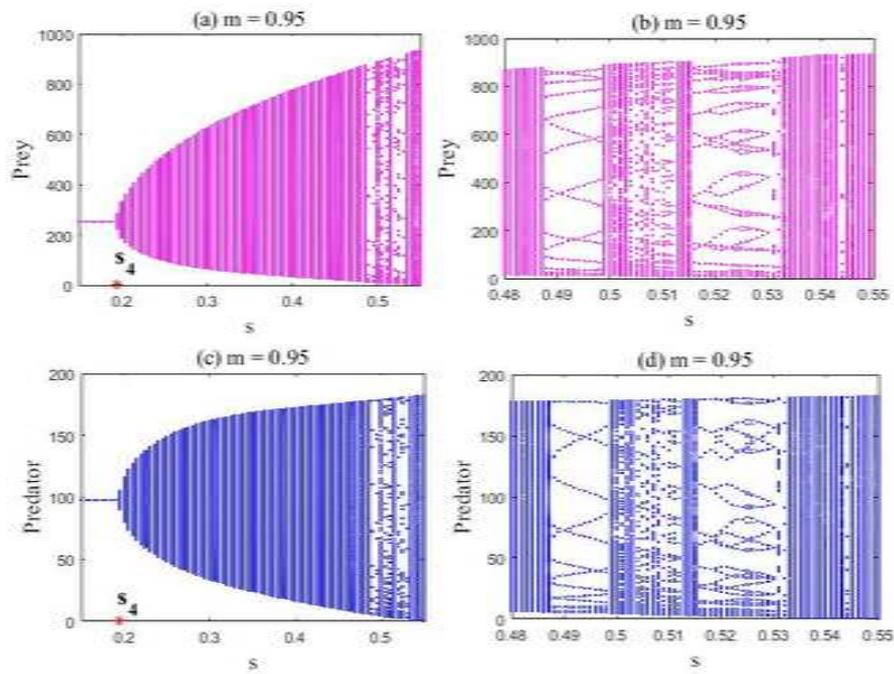}
	\caption{Bifurcation diagrams of fractional-order discrete system (\ref{Discrete model}) with step size $s$ as the bifurcation parameter. Prey and predator populations become unstable as the step size crosses the critical length $s = s_4 = 0.194$ (Figs. (a) and (c)). Figs. (b) and (d) show the local amplification corresponding to (a) and (c), respectively, for $s\in [0.48, 0.55]$. Here $m = 0.95, c=0.45$ and other parameters are as in Fig. 1.}
\end{figure}

\begin{figure}[H]
	\hspace{-1.3cm}
	\includegraphics[width=6in, height=5in]{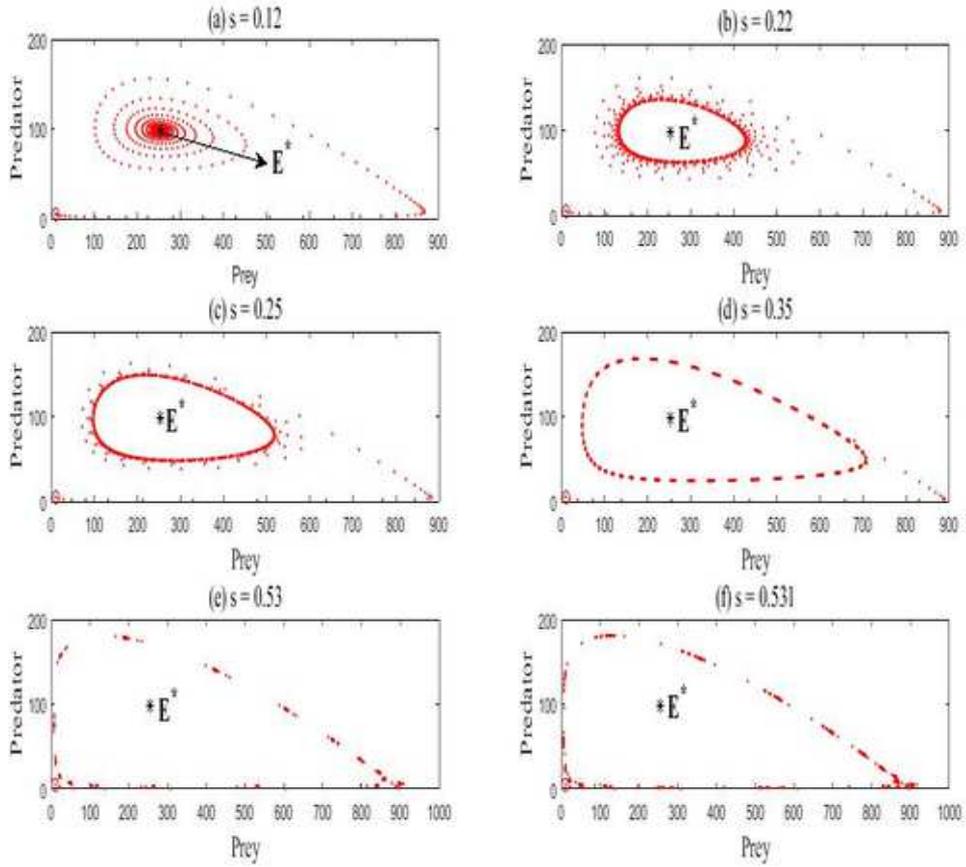}
	\caption{Phase portraits  of fractional-order discrete system (\ref{Discrete model}) for some particular step-size $(s)$ corresponding to Figs. 7(a) and 7(c). Parameters are as in Fig. 7.}
\end{figure}

One can also compute $1+trace(J)+det(J) = 0$ and $\xi_1 = 1 - r \frac{s^m}{m\Gamma (m)} = -1,$ $\xi_2 = 1 + \frac{s^m}{m\Gamma (m)} \bigg(\frac{\theta \alpha (1-c) K}{1+\alpha (1-c)h K} - d\bigg) = 1$ at the predator-free fixed point $E_1$ for $c = c_1$ and $s = s_5 = \sqrt[m]{\frac{2 m\Gamma(m)}{G}} = 0.7279$. Thus, following Remark 3.1, the predator-free fixed point $E_1$ undergoes a flip bifurcation (Fig. 9) at the critical step size $s = s_5$.
\begin{figure}[H]
	\hspace{-2in}
	\includegraphics[width=12in, height=3in]{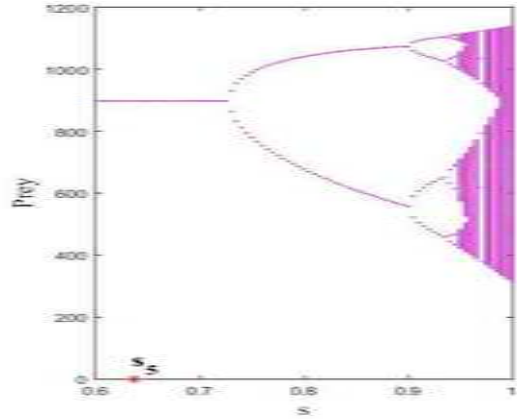}
	\caption{Bifurcation diagram of prey population of fractional-order discrete system (\ref{Discrete model}) with respect to the step size $s$. It shows that the predator-free fixed point $E_1$ containing prey species only undergoes a flip bifurcation at $s = s_5 = 0.7279$. Here $m = 0.95$ and other parameters are as in Fig. 1.}
	\label{Bifurcation.eps}
\end{figure}
\section{Discussion}
This paper generalizes the results of continuous system predator-prey model that considers the effect of habitat complexity. This generalization has been accomplished in two phases. In the first phase, we constructed a fractional-order predator-prey model considering the fractional derivatives in Caputo sense. In the second phase, the fractional-order predator-prey model was discretized. Rigorous mathematical and computational results in relation to the stability of both the systems was presented. Proving existence of Hopf bifurcation with respect to the fractional-order of the derivatives in both discrete and fractional systems is rare in contemporary studies. We have presented it both theoretically and numerically, showing the novelty of this study. For the fractional-order system, we proved different mathematical results like positivity and boundedness, local and global stability of different equilibrium points. It is shown that the trivial equilibrium $E_0$ is always unstable saddle and the predator-free equilibrium is globally asymptotically stable for any value of fractional-order $m \in (0, 1]$ if the degree of habitat complexity exceeds some upper threshold value $c_1$. The solution, however, takes more time to reach the predator-free equilibrium as the value of fraction order is reduced. At the intermediate level of habitat complexity $(c_2<c<c_1)$, the system becomes both locally and globally asymptotically stable around $E^*$ for any value of $m$. These dynamics are consistent with the integer order system $m = 1$. Stability of the interior equilibrium, however, depends on the fractional-order $m$ if the strength of habitat complexity is very low and the system shows order-dependent instability. If $0<c<c_2$ then there exists a critical value $m^*$ of the fractional-order $m \in (0, 1)$ such that the coexistence equilibrium is stable if $m<m^*$ and unstable if $m$ crosses $m^*$. In case of integer order system ($m=1$), the coexistence equilibrium is, however, unstable for all $c\in (0, c_2)$. Simulation results also agree perfectly with the analytical results. Discretization of the fractional-order system was done with piecewise constant arguments and the dynamics of this discrete model was explored. It is observed that the dynamics of the discrete system depends on both the step-size and fractional-order. Existence of Hopf and flip bifurcations have been shown both theoretically and numerically. It is also observed that the discrete fractional-order system shows more complex dynamics as the step size becomes larger. Our simulation results revealed that the discrete system shows period doubling route to chaos for larger step-size.


\bibliography{}   

\end{document}